\documentclass[10pt]{article}

\usepackage{latexsym}
\usepackage{amsfonts}
\usepackage{multicol}
\usepackage{xcolor}

\begin{document}
\parindent=0.2in
\parskip .1in

\color{blue}

\noindent {\Huge {\bf The AMC -- What It Is and Why It Matters}}

\color{black}

\vspace*{.2in}

\noindent {\Large {\bf {\em B\'ela Bajnok}}\footnote{B\'ela Bajnok is a Professor of Mathematics at Gettysburg College and is the Director of the American Mathematics Competitions program of the MAA.  His email address is AMCDirector@MAA.org.}}

\vspace*{.4in}

\begin{multicols}{2}

\noindent 
The first national mathematical competition in the world was organized in Hungary in 1894, several decades before similar events in any other country\footnote{For a comprehensive historical study of national mathematics competitions see {\em Educ. Stud. Math.} 2, 80--114 (1969).}.  At around the same time, the Hungarian journal K\"OMAL\footnote{K\"OMAL is the acronym formed by the first letters of the Hungarian words for high school mathematics competition.} was established, aimed at training high school students in mathematical problem solving, offering year-round competitions, and celebrating student achievements.  Except for the interruptions by the two world wars, both the competition and the journal have been running continuously.

Growing up in Budapest, I was keenly involved in mathematical competitions.  I liked training for them, as I was introduced to intriguing mathematical topics and problems that I didn't see during regular classes; I enjoyed participating in them, even when I didn't perform well; and I loved getting to know many of the other contestants and forming lasting friendships with them.  I have remained involved with mathematical competitions ever since, and thus was excited and honored in 2017, when I became I became the director of the American Mathematics Competitions (AMC) program of the Mathematical Association of America (MAA).

The  AMC consists of a series of mathematical competitions, a collection of curriculum materials, and a variety of engagement opportunities -- all designed to strengthen problem solving skills, foster a love of mathematics, and identify, develop, and nurture talent in middle school and high school students across the United States.  Below I outline these programs in more detail, and describe the roles that they play in the lives of the mathematicians of the future.

The competitions start with the  AMC 8, the  AMC 10, and the  AMC 12 exams, open to students in grade 8 or below, grade 10 or below, and grade 12 or below, respectively.  Over 300,000 students take one of these multiple-choice exams each year.  Afterward, based on their scores on the AMC 10 or 12, about 5000 students are invited to take the  American Invitational Mathematics Exam (AIME), a challenging three-hour exam where the answer to each of the 15 problems is a nonnegative integer under 1000.  The competition series then culminates with the  USA Mathematical Olympiad (USAMO) and the  USA Junior Mathematical Olympiad (USAJMO), offered to approximately 500 students based on their combined performances on the AMC 12 and the AIME or the AMC 10 and the AIME.  These Olympiads follow the style of the  International Mathematics Olympiad (IMO): they ask students to provide rigorous proofs for three problems on each of two consecutive days, with an allowed time of four and a half hours each day.  More information on our competitions can be found at \texttt{https://www.maa.org/math-competitions}.  

Rather than just participating in competitions, we at the AMC believe that students should engage in regular mathematical problem solving activities.  For this reason, the AMC program places an increasing emphasis on providing students and teachers with innovative supporting resources.  Our site  M-Powered at \texttt{www.maa-amc.org}, launched last year, presents a growing supply of materials, such as training videos for educators, developed by MAA's mathematician-at-large James Tanton, that provide insightful explorations of mathematical thinking using problems from AMC competitions.   

We aim to increase the AMC's role as a community builder for students, teachers, friends, and supporters of mathematical problem solving, since opportunities for personal interactions are particularly important for those who feel isolated from others with similar interests.  We have recently introduced events such as AMC Wednesdays and AMC Fireside Chats, and we work closely with the  Art of Problem Solving (AoPS), where participants can engage in discussions on online forums, take AMC competitions at their leisure, and attend live Math Jam sessions with AMC editors.  Our renowned  Math Olympiad Program (MOP) provides intensive in-person training for about 60 students each year, and serves as the pool of potential team members for the USA delegations to the IMO and the  European Girls' Mathematical Olympiad (EGMO)\footnote{In spite of the name, the European Girls' Mathematical Olympiad welcomes teams from all around the world.}.  Our goal is to make these opportunities available to larger and more diverse audiences. 

The importance of building  problem solving skills has been emphasized time and again.  Paul Halmos, for example, is often quoted\footnote{For example, by a report of the Mathematical Association of America Committee on the Teaching of Undergraduate Mathematics (Washington, D.C., 1983) and by the {\em Notices} of the American Mathematical Society (October, 2007, page 1141).} about the importance of problems: 
\begin{quote}  {\sf The major part of every meaningful life is the solution of problems; a considerable part of the professional life of technicians, engineers, scientists, etc., is the solution of mathematical problems.  It is the duty of all teachers, and of teachers of mathematics in particular, to expose their students to problems much more than to facts.}\footnote{``The Heart of Mathematics,'' {\em American Mathematical Monthly} 87 (1980), 519--524.} 
\end{quote}  
Nevertheless, for many students, mathematics still consists mostly of memorizing formulas without understanding their meaning and performing computations without thinking about their purpose.  There are quite a few venues to help students strengthen their problem solving skills, such as math circles, summer camps, and a variety of competitions. The AMC program contributes greatly to this aim in extensive, meaningful, and enjoyable ways.  Bill Serbitz, a math teacher at Wayzata High School in Plymouth, Minnesota, writes:  
\begin{quote} {\sf The AMC provides a source of motivation for my students and my math team members.  It is a nice collection of questions that require a wide variety of skills at one time.  As such, it continues to serve as a great source for developing problem-solving skills at differentiated levels of difficulty.  Because of its accessibility across the globe, It is also useful for building community between like-minded students and teachers worldwide!}
\end{quote}  
Indeed, our goal is to reach every middle and high school student in the US, regardless of background and ability, offering them opportunities to develop their problem solving skills, to challenge themselves, and to learn how to think mathematically.  We hope that the AMC inspires many of them to go on in mathematics and become strong math (or other STEM) majors in college.       

In a radical change from previous practice, two years ago we appointed four large, carefully-chosen  editorial boards that, under the guidance of their editors-in-chief, are in charge of proposing problems, reviewing and rating the problems, assembling the exams, and meticulously editing the problems and their solutions; they also have an opportunity to influence general AMC policies.  This move to a crew of about 120 associate editors allowed us to greatly increase our diversity: now fewer than half of our associate editors are white males.  In addition, we were able to appoint people with different cultural backgrounds, from across the nation as well as from overseas; to include teachers from middle schools and high schools, and also those working in industry; to have a wide range of ages represented, from undergraduates to retired mathematicians; and to have many more fields of the mathematical sciences represented, including statistics.  It is nice to attend editorial board meetings and social events with such a diverse and accomplished group of people.  

The creation of problems suitable for the AMC is a highly challenging task.  
Even at the beginning levels where problems stay close to the standard school curriculum, we aim for problems that ask questions in novel and fun ways.  Coming up with beautiful -- yet still elementary -- problems at the Olympiad level is particularly challenging.  Yet, year after year, collections of ingenious and captivating problems are created; everyone with an interest is encouraged to explore them.\footnote{The full library of previous AMC problems and solutions is available on the Art of Problem Solving website at \texttt{https://artofproblemsolving.com/wiki/index.php/AMC\_Problems\_and\_Solutions}.}    

We live in exceptionally turbulent times, and this has presented the AMC with both  challenges and opportunities.  The COVID-19 pandemic forced us to change how our competitions are administered; while moving our exams online is welcome by many, doing so makes it more difficult or impossible for some students to participate, and raises questions about maintaining the integrity of our competitions.  We also must ask ourselves how to shape our competitions and other programs so they are best suited for our goals.   Should we offer exams that ask students to present arguments, allow them to work in teams, or have longer time-frames -- and if so, in what way?  Do our training programs, such as MOP, align well with our broader objectives?   

Our main challenges, however, are around  diversity and inclusion.  It has long been observed that female students are underrepresented in the AMC program; this underrepresentation becomes increasingly more pronounced at the later rounds of the competition series.  In order to address this troubling issue, we took some recent actions that already appear to be making a difference.  Recognizing that talented female students tend to be involved in more extracurricular activities\footnote{See, Chachra, Chen, Kilgore, and Sheppard, ``Outside the Classroom: Gender Differences in Extracurricular Activities in Engineering Students'', \texttt{https://files.eric.ed.gov/fulltext/ED542117.pdf}.}, we have added math circles and other learning centers as places where our exams can be administered.  We have raised the number of female students invited to MOP, and organize a special pre-MOP program for them.  We have also greatly increased the number of female students receiving awards: due to some recent donations, we will be able to recognize the 20 top-scoring young women with \$1,000 scholarships and another 580 top-scoring women across the country with Certificates of Excellence. During the last two years we also honored the members of our EGMO team alongside the IMO team and the USAMO winners at our annual award ceremony.  Full equality will take time to achieve, but there are signs that we are moving in the right direction; for example, about 14 percent of students taking USAMO 2020 were women, compared to 6--8 percent in the past.  Besides gender equality, we must also focus on other axes of diversity.  We aim to attract more African American and Latinx students by an increased outreach to schools with high percentages of minority students  and by offering the AMC 8, 10, and 12 exams in Spanish.  We also intend to have a more balanced geographic distribution by forming closer connections with MAA sections.

The AMC program would not be possible without the help of many individuals.  I am immensely grateful to our four editorial boards that create our exams each year.  Special thanks go to the AMC Office at the MAA, led by Jennifer Barton,  whose small team is in charge of the operational and administrative aspects of the program.  We are grateful to the brilliant folks at AoPS who provided tailor-made platforms for our competitions and other events to move online.  And we are, of course, profoundly thankful to the thousands of colleagues across the nation who train the students, organize and proctor the exams, and provide us with valuable feedback.             

I would like to close with a request: please consider becoming involved with the AMC, whether working directly with students, helping us with one of our programs or engagement opportunities, or joining one of our editorial boards.  We are in this together; you will be influential not just with the AMC, but with future generations of mathematicians and with mathematics education as a whole.

\end{multicols}

\end{document}